\documentclass[onefignum,onetabnum]{siamart171218}

\ifpdf
\hypersetup{
  pdftitle={Stabilization of the response of cyclically loaded lattice spring models with plasticity},
  pdfauthor={I. Gudoshnikov, O. Makarenkov}
}
\fi

\usepackage{amsfonts}
\usepackage{amssymb}
\usepackage{graphicx}
\usepackage{epstopdf}
\usepackage{algorithmic,comment}
\ifpdf
  \DeclareGraphicsExtensions{.eps,.pdf,.png,.jpg}
\else
  \DeclareGraphicsExtensions{.eps}
\fi

\newsiamremark{remark}{Remark}
\newsiamremark{hypothesis}{Hypothesis}
\crefname{hypothesis}{Hypothesis}{Hypotheses}
\newsiamthm{claim}{Claim}

\headers{Limit cycles of the planar passive biped walking down a slope}{O. Makarenkov}

\title{Existence and stability of a limit cycle in the model of a planar passive biped walking down a slope\thanks{Submitted to the editors DATE.
}}

\author{Oleg Makarenkov\thanks{Deptartment of Mathematical Sciences, University of Texas at Dallas, Richardson, TX 
  (\email{makarenkov@utdallas.edu}, \url{https://www.utdallas.edu/\string~makarenkov/}).}
}

\usepackage{amsopn}

\begin{document}

\sloppy

\maketitle

\begin{abstract}
We consider the simplest model of a passive biped walking down a slope given by the equations of switched coupled pendula (McGeer, 1990). Following the fundamental work by Garcia et al (1998), we view the slope of the ground as a small parameter $\gamma\ge 0$. When $\gamma=0$ the system can be solved in closed form and the existence of a family of limit cycles (i.e. potential walking cycles) can be established explicitly. As observed in Garcia et al (1998), the family of limit cycles disappears when $\gamma$ increases and only isolated asymptotically stable cycles (walking cycles) persist. However, no rigorous proofs of such a bifurcation (often referred to as Melnikov bifurcation) 
have ever been reported. The present paper fills in this gap in the field and offers the required proof. 
\end{abstract}



\section{Introduction} 
\begin{figure}[h]
\begin{center}
\includegraphics[scale=0.5]{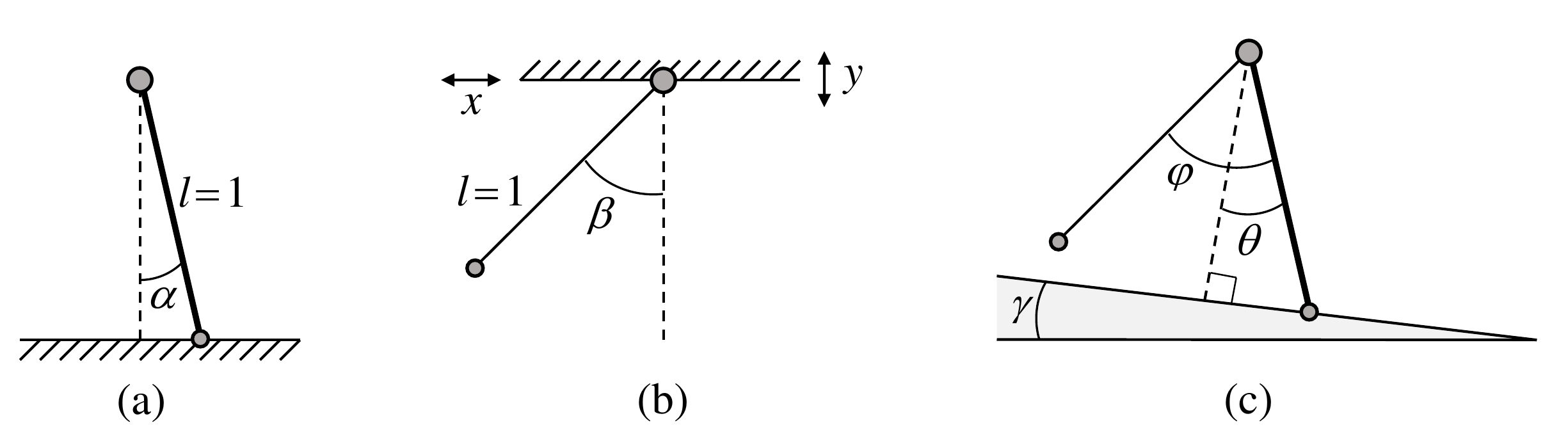}  
\caption{Building blocks ((a) and (b)) of a planar passive walker (c).}\label{walker}
\end{center}
\end{figure}

In his celebrated paper \cite{mcgeer} McGeer proposed to view the passive bipedal walker of Fig.~\ref{walker}a as a combination of a pendulum with a fixed pivot (Fig.~\ref{walker}b)
$$
   \ddot\alpha-g\sin\alpha=0\quad \mbox{(stance leg)}
$$
and a pendulum with a moving pivot (Fig.~\ref{walker}c)
$$
   \ddot \beta+x''\cos\beta+(y''+g)\sin\beta=0 \quad\mbox{(swing leg)},
$$
which gives
\begin{eqnarray}\label{pendula}
\begin{array}{l}
      \ddot \theta-\sin(\theta-\gamma)=0,\\
      \ddot \theta-\ddot\phi+\dot\theta^2\sin\phi-\cos(\theta-\gamma)\sin\phi=0.
      \end{array}
\end{eqnarray}
When the heelstrike occurs (i.e. when $\phi=2\theta$), the stance and swing legs swap their roles and the state vector $(\theta,\dot\theta,\phi,\dot\phi)^T$ jumps as follows 
 \begin{equation}\label{jump}\def\arraystretch{1.2}
 \begin{array}{l} 
\begin{array}{l}
      \left(\hskip-0.07cm\begin{array}{c}
      \theta(t^+)\\
      \dot\theta(t^+)\\
      \phi(t^+)\\
      \dot\phi(t^+)
      \end{array}\hskip-0.07cm\right)=J(\theta(t))\hskip-0.07cm\left(\hskip-0.07cm\begin{array}{c}
      \theta(t^-)\\
      \dot\theta(t^-)\\
      \phi(t^-)\\
      \dot\phi(t^-)
      \end{array}\hskip-0.07cm\right)\hskip-0.07cm, \quad {\rm if}\ \phi(t)=2\theta(t),
   \end{array}
   \end{array}
\end{equation}
where 
$$
J(\theta)=  \left(\hskip-0.07cm\begin{array}{cccc}
      -1 & 0 & 0 & 0\\
      0 & \cos 2\theta & 0 & 0\\
      -2 & 0 & 0 & 0\\
      0 & \begin{array}{c}(1-\cos2\theta)\cos 2\theta \end{array}& 0 & 0
      \end{array}\hskip-0.07cm\right).
$$
Using Newton's method, McGeer found that the switched system (\ref{pendula})-(\ref{jump}) admits a limit cycle, whose period is close to $T=3.8$ for small values of slope $\gamma>0.$ A justification of the existence of such a limit cycle was offered in Garcia el al \cite{garcia}, where the change of the variables 
\begin{equation}\label{change}
   \gamma=\delta^{3/2},\quad \theta(t)=\delta^{1/2}\Theta(t),\quad \phi(t)=\delta^{1/2}\Phi(t)
\end{equation}
is proposed to expand (\ref{pendula})-(\ref{jump}) in the powers of small parameter $\delta>0$ and to investigate the existence of the limit cycle based on the leading order terms. The paper \cite{garcia} offers several important insights linking the properties of the reduced system to the limit cycles of full switched system (\ref{pendula})-(\ref{jump}), but doesn't focus on the rigorous proofs. The goal of the present paper is to provide a rigorous proof of the existence of an attracting limit cycle in (\ref{pendula})-(\ref{jump}) using  appropriate results of the classical perturbation theory. 

\vskip0.2cm

\noindent The paper is organized as follows. In the next section we incorporate the change of the variables (\ref{change}) in switched system (\ref{pendula})-(\ref{jump}) and obtain a switched system (\ref{pendulaeps})-(\ref{jumpeps}) with a small parameter $\delta>0$ (which corresponds to a perturbation term). In Section~\ref{Psection} we follow the idea of Garcia et al \cite{garcia} and introduce a 2-dimensional Poincare map $(\theta,\omega)\mapsto P(\theta,\omega,\delta)$ associated to the perturbed switched system (\ref{pendulaeps})-(\ref{jumpeps}). In Section~\ref{sectiondelta0} we show that, when $\delta=0,$  the Poincare map $(\theta,\omega)\mapsto P(\theta,\omega,\delta)$ admits a family of fixed points $(\theta,\omega)=\xi(s)$, where $\xi\in C^1(\mathbb{R},\mathbb{R}^2)$ and $s$ is a parameter. In this way the problem of the existence of limit cycles to the perturbed switched system (\ref{pendulaeps})-(\ref{jumpeps}) reformulates as a problem of bifurcation of asymptotically stable fixed points to the Poincare map $(\theta,\omega)\mapsto P(\theta,\omega,\delta)$ from the family $(\theta,\omega)=\xi(s)$ as $\delta$ crosses 0. The problem obtained is a classical problem of the theory of nonlinear oscillations coming back to Malkin \cite{malkin} and Melnikov \cite[Ch.~4, \S6]{guck},
 and developed in Loud \cite{loud}, Chicone \cite{chicone}, Rhouma-Chicone \cite{chicone1}, Buica et al \cite{buica1}, Kamenskii et al \cite{kmn}, Makarenkov-Ortega \cite{ortega} and others. In this paper we follow the references \cite{kmn} and \cite{ortega} to provide in Section~\ref{perturbsec} a concise perturbation theorem (Theorem~\ref{theoremmain}) on bifurcation of fixed points from families in Poincare maps. Though the theorem doesn't look new, it seems it has never been formulated in a rigorous form in the literature before. This perturbation theorem is then applied to the Poincare map $(\theta,\omega)\mapsto P(\theta,\omega,\delta)$ of the passive biped in Sections~\ref{stabsection} and \ref{section7}. In Section~\ref{conclusions} (Conclusions) we discuss the value of this work to the field of perturbation theory. The proof of Theorem~\ref{theoremmain} is given in Appendix~\ref{appA} and Appendix~\ref{hfsolution} contains some technical formulas. All symbolic computations have been executed in Wolfram Mathematica 11.3.

\vskip0.2cm

\noindent Despite of extensive literature on bifurcation of fixed points from 1-parameter families, the paper by Glover et al \cite{glover} on large amplitude oscillations in a suspension bridge model seems to be the only example of such a bifurcation accessible to general public. The significant contribution of the present paper is in a rigorous introduction of a one more example of bifurcation from $1$-parameter families that is noticeable to society on the one hand and is well regarded in engineering community on the other hand. 


\section{Expanding McGeer's model of passive biped into the powers of the slope of the ground}

\noindent Incorporating the change of the variables (\ref{change}) into the switched system (\ref{pendula})-(\ref{jump}) and using that
$$
    \sin\tau=\tau-\dfrac{\tau^3}{3!}+\dfrac{\tau^5}{5!}-\dfrac{\tau^7}{7!}+...,\qquad  \cos\tau=1-\dfrac{\tau^2}{2!}+\dfrac{\tau^4}{4!}-\dfrac{\tau^6}{6!}+...,
$$
one gets (see Garcia et al \cite{garcia})
\begin{eqnarray}\label{pendulaeps}
\def\arraystretch{1.8}
&& \hskip0.4cm \begin{array}{l}
   \ddot \Theta-(\Theta-\delta)+\dfrac{1}{6}\delta\Theta^3+o_1(\delta)=0,\\
   \ddot\Theta-\Phi-\ddot\Phi+\delta \dot\Theta^2\Phi+\dfrac{1}{2} \delta\Theta^2\Phi+\dfrac{1}{6}\delta\Phi^3+o_2(\delta)=0,
\end{array}\\ \label{jumpeps}
& & \hskip0.4cm  \begin{array}{l}
      \left(\hskip-0.07cm\begin{array}{c}
      \Theta(t^+)\\
      \dot\Theta(t^+)\\
      \Phi(t^+)\\
      \dot\Phi(t^+)
      \end{array}\hskip-0.07cm\right)=J(\Theta(t),\delta)\hskip-0.07cm\left(\hskip-0.07cm\begin{array}{c}
      \Theta(t^-)\\
      \dot\Theta(t^-)\\
      \Phi(t^-)\\
      \dot\Phi(t^-)
      \end{array}\hskip-0.07cm\right)\hskip-0.07cm,\quad {\rm if}\ \Phi(t)=2\Theta(t),
   \end{array}
\end{eqnarray}
where 
$$
J(\Theta,\delta)=\left(\hskip-0.07cm\begin{array}{cccc}
      -1 & 0 & 0 & 0\\
      0 & \begin{array}{c}
      1-\dfrac{1}{2}\delta(2\Theta)^2 +o_3(\delta)
      \end{array}& 0 & 0\\
      -2 & 0 & 0 & 0\\
      0 & \left(1-\dfrac{1}{2}\delta(2\Theta)^2 +o_3(\delta)\right)\left(\dfrac{1}{2}\delta (2\Theta)^2+o_4(\delta)\right)& 0 & 0
      \end{array}\hskip-0.07cm\right)
$$
and $o_i(\delta)$ stay for the remainders (perhaps dependent on $\Theta$ and $\Phi$) such that $o_i(\delta)/\delta\to 0$ as $\delta\to 0$ uniformly with respect to $(\Theta,\Phi)$ from any compact set. 

\section{\boldmath{The Poincare map $(\theta,\omega)\mapsto P(\theta,\omega,\delta)$ induced by the heelstrike threshold}} \label{Psection} To construct the Poincare map induced by the hyperplane $\Phi=2\Theta,$ we will consider the initial condition 
$(\Theta(t^+),
      \dot\Theta(t^+),
      \Phi(t^+),
      \dot\Phi(t^+))^T$ given by (\ref{jumpeps}). Because of the properties of the matrix $J(\Theta,\delta)$ any vector $(\Theta(t^+),
      \dot\Theta(t^+),
      \Phi(t^+),
      \dot\Phi(t^+))^T$ coming from (\ref{jumpeps}) has the form
\begin{equation}\label{reductionform}
   (\Theta(t^+),
      \dot\Theta(t^+),
      \Phi(t^+),
      \dot\Phi(t^+))=\left(\theta,\omega,2\theta,\left(2\delta \theta^2+o_4(\delta)\right)\omega\right).
\end{equation}
In other words, knowing  that $(\Theta(t^+),\dot\Theta(t^+))=(\theta,\omega)$, we can use formula (\ref{reductionform}) to obtain the respective values of 
      $\Phi(t^+)$ and 
      $\dot\Phi(t^+)$. Defining 
$$
   \Delta(\theta,\omega,\delta)=\left(\begin{array}{c}
   -\theta\\
   (1-2\delta \theta^2+o_3(\delta))\omega\end{array}\right),
$$      
      we can introduce a 2-dimensional Poincare map as follows
\begin{equation}\label{P}
   P(\theta,\omega,\delta)=\Delta\left[
   \left(\Theta, 
   \dot\Theta\right)(T(\theta,\omega,\delta),\theta,\omega,\delta),\delta\right],
\end{equation}
where 
$t\mapsto(\Theta,\dot\Theta,\Phi,\dot\Phi)^T(t,\theta,\omega,\delta)$ is the solution of 
(\ref{pendulaeps}) with the initial condition 
\begin{equation}\label{initialcond}
(\Theta(0),
      \dot\Theta(0),
      \Phi(0),
      \dot\Phi(0))=\left(\theta,\omega,2\theta,2\delta\theta^2\omega+o_4(\delta)\omega\right)
\end{equation} 
and $T(\theta,\omega,\delta)$ is the time satisfying 
      \begin{equation}\label{TT}
      \begin{array}{l}
         \Phi(T(\theta,\omega,\delta),\theta,\omega,\delta)=2\Theta(T(\theta,\omega,\delta),\theta,\omega,\delta),\\
         \Phi(t,\theta,\omega,\delta)\not=2\Theta(t,\theta,\omega,\delta),\quad t\in(0,T(\theta,\omega,\delta)).
         \end{array}
      \end{equation}
      
\section{\boldmath{Families of fixed points of the Poincare map $(\theta,\omega)\mapsto P(\theta,\omega,\delta)$ when $\delta=0$}}\label{sectiondelta0}
When $\delta=0$, the system (\ref{pendulaeps}) and the initial condition (\ref{reductionform}) take the form
\begin{equation}\label{pendulareduced}
 \hskip0.4cm \begin{array}{l}
   \ddot \Theta-\Theta=0,\\
   \ddot\Theta-\Phi-\ddot\Phi=0,
\end{array}\qquad (\Theta(0),\dot\Theta(0),\Phi(0),\dot\Phi(0))=(\theta,\omega,2\theta,0),
\end{equation}
whose solution is
\begin{equation}\label{f1}
\begin{array}{ll}
\Theta(t,\theta,\omega,0)= \dfrac{1}{2}e^{-t}(1+e^{2t})\theta+\dfrac{1}{2}e^{-t}(-1+e^{2t})\omega,\\
\Phi(t,\theta,\omega,0)=2\theta\cos t+\dfrac{e^{-t}}{4}(1+e^{2t}-2e^t\cos t)\theta+\dfrac{e^{-t}}{4}(e^{2t}-1-2e^t\sin t)\omega.
\end{array}
\end{equation}
Observe that $(\theta,\omega)=P(\theta,\omega,0)$ if and only if
\begin{equation}\label{ifif}
\begin{array}{l}
   \theta=-\Theta(T,\theta,\omega,0),\\
   \omega=\dot\Theta(T,\theta,\omega,0),\\
   \Phi(T,\theta,\omega,0)=2\Theta(T,\theta,\omega,0).
\end{array}
\end{equation}
The first two equations of (\ref{ifif}) give
\begin{equation}\label{omegatheta}
    \omega=\alpha(T)\theta,\quad{\rm where}\quad \alpha(T)=-\dfrac{1+e^T}{-1+e^T}.
\end{equation}
Substituting (\ref{omegatheta}) into the third equation of (\ref{ifif}) one obtains the following equation for $T$
\begin{equation}\label{T}
   -3+3e^T+3(-1+e^T)\cos T+\sin T+e^T\sin T=0,
\end{equation}
whose roots on $(0,2\pi)$ 
\begin{equation}\label{f2}
   T_1=\pi,\quad T_2=3.81209... .
\end{equation}
According to Garcia et al \cite{garcia} only roots within the $(0,2\pi)$ correspond to ``reasonably anthropomorphic gaits''. Also, following   
Garcia et al \cite{garcia}, we will stick to the second root $T_2$ because it corresponds to a symmetric gait in the following sense: plugging $\omega=\alpha(T_2)\theta$ into the third equation of (\ref{ifif}) gives approximately
$$
   -1.5339 e^{-t}+0.0339021 e^t+1.5\cos t+0.522601\sin t=0,
$$
whose only solution on $(0,T_2)$ is $T_2/2$ where one has  \begin{equation}\label{slip}
\Theta(T_2/2,\theta,\omega,0)=\Phi(T_2/2,\theta,\omega,0)=0.
\end{equation}
The property (\ref{slip}) corresponds to the event where the two legs coincide. Though (\ref{slip}) formally implies a heel-strike (the third equation of (\ref{ifif}) holds at $T=T_2/2$), it corresponds to just grazing of the swing leg through the floor and no impact event physically occurs. If the value of $\gamma$ increases, then, formally speaking, an impact occurs at $T=T_2/2$, but we will still ignore the impact coming from $T=T_2/2$ as motivated by the experiments (the actual experimental passive planar walker makes slight swings in the 3rd dimension which rules out the impact at $T=T_2/2$, see \cite{video}). In other words, for the reasons just explained and following Garcia et al \cite{garcia}, we will consider the Poincare map (\ref{P}) with 
$$
  T(\theta,\omega,\delta)\to T_2\quad{\rm as}\quad\delta\to 0
$$
which satisfies the first condition of (\ref{TT}) even though it ``slightly'' violates the second condition of (\ref{TT}) in the neighborhood of $T_2/2.$



\section{Perturbation theorem for two-dimensional Poincare maps} \label{perturbsec}
\noindent Throughout this section we assume that the unperturbed Poincare map $(\theta,\omega)\mapsto P(\theta,\omega,0)$ admits 
a family of fixed points, i.e. $P(\xi(s),0)=\xi(s)$ for all $s\in\mathbb{R}$, where $s\mapsto \xi(s)$ is a $C^1$ curve. 
Note, the latter property implies that $P_{(\theta,\omega)}(\xi(s),0)\xi'(s)=\xi'(s),$ which means that 
one of the eigenvalues of the matrix $P_{(\theta,\omega)}(\xi(s),0)$ is always 1 for all $s\in\mathbb{R}.$ To make the notations less bulky we will  identify $P(\theta,\omega,\delta)$ with $P((\theta,\omega),\delta)$ as it doesn't seem to cause any confusion.

\vskip0.2cm

\noindent Fix some $s_0\in\mathbb{R}$ and put
$$
   (\theta_0,\omega_0)=\xi(s_0).
$$
Denote by 
$y$ and $\tilde y$ the eigenvectors of 
$P_{(\theta,\omega)}(\theta_0,\omega_0,0)$ that correspond to the eigenvalues 1 and 
$\rho\not=1
$ respectively. We then  denote by $z$ and $\tilde z$ the eigenvalues of  $P_{(\theta,\omega)}(\theta_0,\omega_0,0)^T$ that correspond to the eigenvalues $1$ and $\rho\not=1,$ and such that
\begin{equation}\label{prop1}
   z^Ty=\tilde z^T\tilde y=1.
\end{equation}
It can be verified that 
\begin{equation}\label{prop0}
    z^T\tilde y=\tilde z^T y=0.
\end{equation}
Properties (\ref{prop1}) and (\ref{prop0}) imply that
\begin{equation}\label{decomposition}
   \zeta=z^T\zeta y+\tilde z^T\zeta\tilde y,\quad\mbox{for any}\ \zeta\in\mathbb{R}^2.
\end{equation}
We will also assume that $z$ doesn't depend on the choice of $s_0$, in which case we have
\begin{equation}\label{dd}
    z^T(P_{(\theta,\omega)}(\xi(s),0)-I)=0,\quad\mbox{for all}\ s\in\mathbb{R}.
\end{equation}

\noindent The following theorem is a corollary of the results of  Kamenski et al \cite{kmn} and  Makarenkov-Ortega \cite{ortega}.  

\begin{theorem} \label{theoremmain} Let $P$ be a $C^3$ function. If, for each $\delta\in\mathbb{R}$, the Poincare map
$(\theta,\omega)\mapsto P(\theta,\omega,\delta)$ admits
a fixed point $(\theta_\delta,\omega_\delta)$ such that 
\begin{equation}\label{convergence}
  (\theta_\delta,\omega_\delta)\to (\theta_0,\omega_0)\quad {\rm as}\quad \delta\to 0,
\end{equation}
then
\begin{equation}\label{necessary}
   z^T P_\delta(\theta_0,\omega_0,0)=0.
\end{equation}
Assume that the eigenvector $z$ of 
$P_{(\theta,\omega)}(\theta_0,\omega_0,0)^T$ that corresponds to the eigenvalue 1 doesn't depend on $s_0.$ 
If, in addition to (\ref{necessary}),  it holds that
\begin{equation}\label{sufficient}
   z^T (P_\delta)_{(\theta,\omega)}(\theta_0,\omega_0,0)y\not=0,
\end{equation}
then, for all $|\delta|$ sufficiently small, the Poincare map  
$(\theta,\omega)\mapsto P(\theta,\omega,\delta)$ does indeed have a fixed point $(\theta_\delta,\omega_\delta)$ that satisfies (\ref{convergence}). 
The fixed point $(\theta_\delta,\omega_\delta)$ is asymptotically stable, if the eigenvalue 
$\rho\not=1$ of $P_{(\theta,\omega)}(\theta_0,\omega_0,0)$ satisfies
\begin{equation}\label{stab1}
|\rho|<1,
\end{equation}
and if 
(\ref{sufficient}) holds in the stronger sense
\begin{equation}\label{stab2}
   z^T (P_\delta)_{(\theta,\omega)}(\theta_0,\omega_0,0)y<0.
\end{equation}
\end{theorem}

\section{\boldmath Stability of the family 
$\omega=\alpha(T_2)\theta$
of fixed points of the Poincare map $(\theta,\omega)\mapsto P(\theta,\omega,\delta)$ corresponding to $\delta=0$}\label{stabsection}

\noindent As explained in Section~\ref{perturbsec}, one of the eigenvalues of  matrix $P_{(\theta,\omega)}(\theta,\alpha(T_2)\theta,0)$ is always 1. In this section we compute the second eigenvalue (named $\rho$) of $P_{(\theta,\omega)}(\theta,\alpha(T_2)\theta,0)$ and verify condition (\ref{stab1}) of Theorem~\ref{theoremmain}. We will see that $\rho$ doesn't depend on $\theta$, so we write $\rho$ as opposed to $\rho(\theta)$ from the beginning.

\vskip0.2cm

\noindent Differentiating  (\ref{P}) with respect to the vector variable $(\theta,\omega),$
$$
\begin{array}{rcl}
P_{(\theta,\omega)}(\theta,\alpha(T_2)\theta,0)&=&\Delta_0 \left(\hskip-0.15cm\begin{array}{c}
\Theta\\
\dot\Theta
\end{array}\hskip-0.15cm\right)_t(T_2,\theta,\alpha(T_2)\theta,0)\hskip0.1cmT_{(\theta,\omega)}(\theta,\alpha(T_2)\theta,0)+\\
&&\qquad\qquad+\Delta_0\left(\hskip-0.15cm\begin{array}{c}
\Theta\\
\dot\Theta
\end{array}\hskip-0.15cm\right)_{(\theta,\omega)}(T_2,\theta,\alpha(T_2)\theta,0),
\end{array}
$$
where
$$
\Delta_0=\left(\hskip-0.07cm\begin{array}{cc}
      -1 & 0 \\
      0 & 1
            \end{array}\hskip-0.07cm\right).
$$
Using formulas (\ref{f1}) and (\ref{f2}) one gets
$$
\left(\hskip-0.15cm\begin{array}{c}
\Theta\\
\dot\Theta
\end{array}\hskip-0.15cm\right)_t(T_2,\theta,\omega,0)=\left(\hskip-0.15cm\begin{array}{c}
\Theta_t(T_2,\theta,\omega,0)\\
\Theta_{tt}(T_2,\theta,\omega,0)
\end{array}\hskip-0.15cm\right)=
\left(\begin{array}{c}
22.6114 \omega + 22.6335 \theta\\
22.6335 \omega + 22.6114 \theta
\end{array}\right)
$$
and so
$$
\left(\hskip-0.15cm\begin{array}{c}
\Theta\\
\dot\Theta
\end{array}\hskip-0.15cm\right)_t(T_2,\theta,\alpha(T_2)\theta)=\theta\left(\begin{array}{c}
-1.0452 \\
-1
\end{array}\right).
$$
In the same way, 
$$
\left(\hskip-0.15cm\begin{array}{c}
\Theta\\
\dot\Theta
\end{array}\hskip-0.15cm\right)_{(\theta,\omega)}(\tau)=\left(\hskip-0.15cm\begin{array}{cc}
\Theta_\theta(\tau) & \Theta_\omega(\tau)\\
\Theta_{t\theta}(\tau) &  \Theta_{t\omega}(\tau)
\end{array}\hskip-0.15cm\right)=\left(\hskip-0.15cm\begin{array}{cc}
22.6335& 22.6114\\22.6114 & 22.6335
\end{array}\hskip-0.15cm\right),
$$
where a shortcut 
$$\tau=(T_2,\theta,\alpha(T_2)\theta,0)$$
is used.
The formula for the derivative of the implicit function (see e.g. Zorich \cite[Sec. 8.5.4 Theorem 1]{zorich}) further yields
\begin{equation}\label{f3}
  T_{(\theta,\omega)}(\theta,\alpha(T_2)\theta,0)=-\left(F_t(T_2,\theta,\alpha(T_2)\theta)\right)^{-1}F_{(\theta,\omega)}(T_2,\theta,\alpha(T_2)\theta),
\end{equation}
where
\begin{equation}\label{formulaF}
  F(t,\theta,\omega)=\Phi(t,\theta,\omega,0)-2\Theta(t,\theta,\omega,0).
\end{equation}
Plugging formulas (\ref{f1}) and (\ref{f2}) into (\ref{f3}), the function $T_{(\theta,\omega)}(\theta,\alpha(T_2)\theta,0)$ computes as
$$
T_{(\theta,\omega)}(\theta,\alpha(T_2)\theta,0)=
\dfrac{1}{\theta}(16.8032,16.0765).
$$
Combining the above findings together we finally get
\begin{equation}\label{Pthetaomega}
P_{(\theta,\omega)}(\theta,\alpha(T_2)\theta,0)=
\left(\begin{array}{cc}
 -5.07075 & -5.8082\\
 5.8082 & 6.55701
 \end{array}\right)
\end{equation}
whose eigenvalues are $1$ and 
$$
\rho=0.48626,
$$ 
so that condition (\ref{stab1}) holds.

\section{\boldmath Bifurcation of isolated fixed points of the Poincare map $(\theta,\omega)\mapsto P(\theta,\omega,\delta)$ from the family $\omega=\alpha(T_2)\theta$ when $\delta$ crosses $0$} \label{section7}

In this section we verify the remaining conditions (\ref{necessary}), (\ref{sufficient}) and (\ref{stab2}) of Theorem~\ref{theoremmain}.

\subsection{\boldmath{Computing $P_\delta$}}
\label{compPdelta}

\noindent Differentiating  (\ref{P}) with respect to $\delta,$ one gets
 \begin{equation}\label{Pdelta}
 \begin{array}{rcl}
   P_\delta(\theta,\omega,0)&=&\Delta_0\hskip-0.05cm\left(\hskip-0.15cm\begin{array}{c}
   \Theta \\
   \dot\Theta\end{array}\hskip-0.15cm\right)_{\hskip-0.05cm t}\hskip-0.05cm(T(\theta,\omega,0),\theta,\omega,0)\hskip0.05cm T_\delta(\theta,\omega,0)+\\
   &&\qquad+\Delta_0\hskip-0.05cm\left(\hskip-0.15cm\begin{array}{c}
   \Theta \\
   \dot\Theta\end{array}\hskip-0.15cm\right)_{\hskip-0.05cm\delta}\hskip-0.05cm(T(\theta,\omega,0),\theta,\omega,0)+\\
   &&\qquad\qquad\qquad+\Delta_\delta\left(
   \left(\Theta, 
   \dot\Theta\right)(T(\theta,\omega,0),\theta,\omega,0),0\right).
\end{array}
\end{equation}

\noindent The terms $\Delta_0$ and $\left(\hskip-0.15cm\begin{array}{c}
   \Theta \\
   \dot\Theta\end{array}\hskip-0.15cm\right)_{\hskip-0.05cm t}\hskip-0.05cm(\tau)$ were computed in the previous section. For the terms $\Delta_\delta(\theta,\omega,0)$ and $\left(\hskip-0.15cm\begin{array}{c}
   \Theta \\
   \dot\Theta\end{array}\hskip-0.15cm\right)\hskip-0.05cm(\tau)$, the definition of $\Delta(\theta,\omega,\delta)$ and formula (\ref{ifif}) yield 
   $$
   \Delta_\delta(\theta,\omega,0)=\left(\begin{array}{c}
   0 \\
   -2\theta^2\omega
\end{array}\right),\qquad
     \left(\hskip-0.15cm\begin{array}{c}
   \Theta \\
   \dot\Theta\end{array}\hskip-0.15cm\right)\hskip-0.05cm(\tau)=\left(\begin{array}{c}
   -\theta \\
   \alpha(T_2)\theta\end{array}\right).
   $$
To compute $T_\delta(\theta,\omega,0)$ we can use function $F$ of the previous section, which gives
\begin{equation}\label{Tdelta}
  T_\delta(\theta,\omega,0)=-\left(F_t(T_2,\theta,\omega)\right)^{-1}F_\delta(T_2,\theta,\omega).
\end{equation}
So it remains to compute the function $t\mapsto\left((\Theta,\Phi)^T\right)_\delta(t,\theta,\omega,0),$ which can be found as the solution $t\mapsto (h(t),f(t))^T$ of the $\delta$-derivative of the initial-value problem (\ref{pendulaeps}) and (\ref{initialcond}):
\begin{equation}\label{linearization} \def\arraystretch{1.5}
\begin{array}{l}
    \ddot h-h+1+\dfrac{1}{6}\Theta(t,\sigma)=0,\\
    \ddot h-f-\ddot f+\left(\dot \Theta(t,\sigma)\right)^2\Phi(t,\sigma)+\dfrac{1}{2}\left(\dot \Theta(t,\sigma)\right)^2\Phi(t,\sigma)+\dfrac{1}{6}\left(\Phi(t,\sigma)\right)^3=0,\\
    h(0)=0,\ \ \dot h(0)=0,\ \ f(0)=0,\ \ \dot f(0)=2\theta^2\omega,
\end{array}
\end{equation}
where $\sigma$ is a shortcut for $\sigma=(\theta,\omega,0).$ After plugging (\ref{f1}) into (\ref{linearization}) we get a system of linear inhomogeneous differential equations, whose solution 
$t\mapsto (h(t),f(t))^T$ is given in Appendix~\ref{hfsolution}. In particular, plugging $t=T_2,$ one gets
\begin{equation}\label{byanalogy}
\begin{array}{l}
 \left(\begin{array}{c}\Theta\\ \Phi\end{array}\right)_\delta(T_2,\theta,\omega,0)=\left(\begin{array}{c}
  h(T_2)\\ f(T_2)\end{array}\right)=\\
   =\left(\begin{array}{c}
  -21.6335 - 236.869 \omega^3 - 717.864 \omega^2 \theta - 
 726.524 \omega \theta^2 - 246.471 \theta^3
\\
 -11.7085 + 669.091 \omega^3 + 1793.6 \omega^2 \theta + 
 1582.73 \omega \theta^2 + 458.155 \theta^3
  \end{array}\right).
\end{array}
\end{equation}

\noindent and
$$
  \left(\begin{array}{c}\Theta\\ \Phi\end{array}\right)_\delta(\tau)=\left(\begin{array}{c}
  h(T_2)\\ f(T_2)\end{array}\right)=\left(\begin{array}{c}
  -21.6335 - 0.871197 \theta^3\\
  -11.7085 - 0.697524 \theta^3
  \end{array}\right).
$$
Formula (\ref{Tdelta}) then provides
$$
   T_\delta(\theta,\omega,0)=\dfrac{0.940403 + 34.0548 \omega^3 + 96.2296 \omega^2 \theta + 
 90.4622 \omega \theta^2 + 28.3414 \theta^3}{
\omega + 0.982912 \theta}.
$$
Plugging all the above findings into formula (\ref{Pdelta}), we  conclude
$$
  P_\delta(\theta,\alpha(T_2)\theta,0)=\left(\begin{array}{c}
  5.85426 + 0.348762 \theta^3\\
  -7.51458 + 1.75673 \theta^3
  \end{array}\right).
$$
 
\subsection{\boldmath{Computing $(\theta_0,\omega_0)$ that satisfies the necessary condition (\ref{necessary})}}
\label{compnecessary} Computing an eigenvector $z$ of the transpose of the matrix (\ref{Pthetaomega}) for the eigenvalue 1, we get 
$$
  z=(-0.69131,-0.722559)^T.
$$
 Therefore, taking into account the relation (\ref{f2}) between $\theta_0$  and $\omega_0$,  the necessary condition (\ref{necessary}) takes the form
 $$
    1.38262-1.51044(\theta_0)^3=0.
 $$  
The solution of this equation is 
 $$
     \theta_0=0.970956,
 $$
 which coincides with the finding of Garcia et al \cite{garcia} (see the table at \cite[p.~15]{garcia}).

\subsection{\boldmath{Computing $P_{\delta(\theta,\omega)}$}} \label{compPdeltathetaomega}

Differentiating (\ref{Pdelta}) with respect to $(\theta,\omega),$ one gets
$$
 \hskip-0.05cm\begin{array}{l}
   P_{\delta(\theta,\omega)}(\theta,\alpha(T_2)\theta,0)\hskip-0.02cm=\hskip-0.02cm\Delta_0\hskip-0.1cm\left[\left(\hskip-0.15cm\begin{array}{c}
   \Theta \\
   \dot\Theta\end{array}\hskip-0.15cm\right)_{\hskip-0.05cm tt}\hskip-0.15cm(\tau)\hskip0.05cm T_{(\theta,\omega)}(\theta,\alpha(T_2)\theta,0)+\right.\\
   \hskip5cm\left.+\hskip-0.05cm\left(\hskip-0.15cm\begin{array}{c}
   \Theta \\
   \dot\Theta\end{array}\hskip-0.15cm\right)_{\hskip-0.05cm t(\theta,\omega)}\hskip-0.15cm(\tau)\hskip0.05cm\right]\hskip-0.1cm T_\delta(\theta,\alpha(T_2)\theta,0)+\\
 \hskip3.15cm+\Delta_0\hskip-0.05cm\left(\hskip-0.15cm\begin{array}{c}
   \Theta \\
   \dot\Theta\end{array}\hskip-0.15cm\right)_{\hskip-0.05cm t}\hskip-0.05cm(\tau)\hskip0.05cm T_{\delta(\theta,\omega)}(\theta,\alpha(T_2)\theta,0)+\\
   \hskip3.15cm+\Delta_0\hskip-0.1cm\left[\left(\hskip-0.15cm\begin{array}{c}
   \Theta \\
   \dot\Theta\end{array}\hskip-0.15cm\right)_{\hskip-0.05cm \delta t}\hskip-0.15cm(\tau)T_{(\theta,\omega)}(\theta,\alpha(T_2)\theta,0)\hskip-0.05cm+\hskip-0.05cm\left(\hskip-0.15cm\begin{array}{c}
   \Theta \\
   \dot\Theta\end{array}\hskip-0.15cm\right)_{\hskip-0.05cm \delta(\theta,\omega)}\hskip-0.15cm(\tau)\hskip0.05cm\right]
   +\\
   \hskip3.15cm+\Delta_{\delta(\theta,\omega)}\left(
   \left(\Theta, 
   \dot\Theta\right)(T(\theta,\alpha(T_2)\theta,0),\theta,\omega,0),0\right)\circ\\
   \hskip5cm\circ\left[\left(\hskip-0.15cm\begin{array}{c}
   \Theta \\
   \dot\Theta\end{array}\hskip-0.15cm\right)_{\hskip-0.05cm t}\hskip-0.15cm(\tau)T_{(\theta,\omega)}(\theta,\alpha(T_2)\theta,0)\hskip-0.05cm+\hskip-0.05cm\left(\hskip-0.15cm\begin{array}{c}
   \Theta \\
   \dot\Theta\end{array}\hskip-0.15cm\right)_{\hskip-0.05cm (\theta,\omega)}\hskip-0.15cm(\tau)\hskip0.05cm\right].
\end{array}
$$
 The terms $\left(\hskip-0.15cm\begin{array}{c}
   \Theta \\
   \dot\Theta\end{array}\hskip-0.15cm\right)_{tt}(t,\theta,\omega,0)$ and $\left(\hskip-0.15cm\begin{array}{c}
   \Theta \\
   \dot\Theta\end{array}\hskip-0.15cm\right)_{t(\theta,\omega)}(t,\theta,\omega,0)$ come by taking the derivatives of $\left(\hskip-0.15cm\begin{array}{c}
   \Theta \\
   \dot\Theta\end{array}\hskip-0.15cm\right)_{\hskip-0.05cm t}\hskip-0.05cm(t,\theta,\omega,0)$ with respect to $t$ and $(\theta,\omega).$ The formulas for  $T_{(\theta,\omega)}(\theta,\alpha(T_2)\theta,0)$ and $T_\delta(\theta,\omega,0)$ were computed in Sections~\ref{stabsection} and \ref{compPdelta}. To compute $T_{\delta(\theta,\omega)}$ we just differentiate the formula for $T_\delta(\theta,\omega,0)$ of Section~\ref{compPdelta} with respect to $(\theta,\omega)$ obtaining
$$
\begin{array}{l}
  \hskip-0.1cm T_{\delta(\theta,\omega)}(\theta,\omega,0)=\left(\dfrac{-0.924333 + 62.7568 \omega^3 + 180.924 \omega^2 \theta + 
 173.941 \omega \theta^2 + 
 55.7142 \theta^3}{(\omega + 0.982912 \theta)^2},\right.\hskip-0.5cm\\
\hskip2.7cm \left.\dfrac{-0.940403 + 68.1095 \omega^3 + 196.648 \omega^2 \theta + 
 189.17 \omega \theta^2 + 
 60.575 \theta^3}{(\omega + 0.982912 \theta)^2}\right).
 \end{array}
$$
By analogy with (\ref{byanalogy}) we compute
$$
\begin{array}{l}
 \left(\begin{array}{c}\Theta\\ \dot\Theta\end{array}\right)_{\delta t}(T_2,\theta,\omega,0)=\left(\begin{array}{c}
  \dot h(T_2)\\ \ddot h(T_2)\end{array}\right)=\\
   \qquad \quad =\left(\begin{array}{c}
  -22.6114 - 717.864 \omega^3 - 2163.65 \omega^2 \theta - 
 2175.14 \omega \theta^2 - 730.293 \theta^3
\\
 -22.6335 - 2163.65 \omega^3 - 6503.87 \omega^2 \theta - 
 6518.19 \omega \theta^2 - 2178.91 \theta^3
  \end{array}\right).
\end{array}
$$   
It remains to find $\Delta_{\delta(\theta,\omega)}(\theta,\omega,0)$ which computes as
$$
  \Delta_{\delta(\theta,\omega)}(\theta,\omega,0)=\left(\begin{array}{cc}
    0 & 0 \\
    -4\theta\omega & -2\theta^2\end{array}\right).
$$  
Combining all the findings together, the matrix $P_{\delta(\theta,\omega)}(\theta,\alpha(T_2)\theta,0)$ finally computes as
$$
   P_{\delta(\theta,\omega)}(\theta,\alpha(T_2)\theta,0)=\dfrac{1}{\theta}\left(\begin{array}{cc}
   218.645 - 3.99563 \theta^3 & 
 209.189 - 4.82387 \theta^3\\
 -21.9105 - 
  94.324 \theta^3 & -20.9629 - 95.2869 \theta^3\end{array}\right).
$$
   
\subsection{Verifying the stability condition (\ref{stab2})} To verify condition (\ref{stab2}), it remains to compute the eigenvector $y$ matrix (\ref{Pthetaomega}) which corresponds to the eigenvalue 1 and satisfies the normalization property (\ref{prop1}) with the vector $z$ of Section~\ref{compnecessary}. Such a computation leads to
$$
y=(15.6468, -16.3541)^T.
$$
Using the formula for $P_{\delta(\theta,\omega)}(\theta_0,\alpha(T_2)\theta_0,0)$ of Section~\ref{compPdeltathetaomega} and the value $\theta_0$ given by Section~\ref{compnecessary}, we get
$$
z^TP_{\delta(\theta,\omega)}(\theta_0,\alpha(T_2)\theta_0,0)y=-2.95323,
$$
so that both the conditions (\ref{sufficient}) and (\ref{stab2}) hold.

\section{Conclusions} \label{conclusions}  In this paper we built upon the fundamental paper by Garcia et al \cite{garcia} and then used the results by Kamenskii et al \cite{kmn} and Makarenkov-Ortega \cite{ortega} in order to offer a step-by-step guide as for how the classical perturbation theory needs to be applied in order to establish the existence and stability of a walking limit cycle in a model of passive biped by McGeer \cite{mcgeer}. Since the dynamics of a passive walker constitutes an important building block of more complex robotics models (engineers use the passive walker dynamics to diminish the energy required for locomotion), we like to think that the present work will stimulate the use of perturbation theory in the field of robotics.

   
\appendix\section
{Derivation of the perturbation theorem of Section~\ref{perturbsec} from the results of Kamenskii et al \cite{kmn} and Makarenkov-Ortega \cite{ortega}} \label{appA}

\noindent The following two results have been established in Kamenskii et al \cite{kmn} and they will play the central role in the perturbation theorem (Theorem~\ref{theoremmain}) that this section develops. We now reformulate the required results of \cite{kmn} in the notations of the present paper to avoid confusion. 

\begin{theorem} \label{theoremkmn} {\rm (two-dimensional version of a combination of \cite[Theorem 1]{kmn} and \cite[Remark 2]{kmn})} Consider a $C^2$-function  $(\theta,\omega,\delta)\mapsto F(\theta,\omega,\delta)$. 
Let $\Pi:\mathbb{R}^2\to\mathbb{R}^2$ be a linear projector invariant with respect to $F_{(\theta,\omega)}(\theta_0,\omega_0,0)$ with $F_{(\theta,\omega)}(\theta_0,\omega_0,0)$ invertible on $(I-\Pi)\mathbb{R}^2.$ Assume that $\Pi F_\delta(\theta_0,\omega_0,0)=0$, $\Pi F_{(\theta,\omega)^2}(\theta_0,\omega_0,0)\Pi h_1\Pi h_2=0$ for any $h_1,h_2\in\mathbb{R}^2$, and that
\begin{equation}\label{invertible}
\begin{array}{l}
     -\Pi F_{(\theta,\omega)^2}(\theta_0,\omega_0,0)h+\Pi (F_\delta)_{(\theta,\omega)}(\theta_0,\omega_0,0),\\
    \qquad \qquad {\rm where}\ h=(I-\Pi)\left(\left.F_{(\theta,\omega)}(\theta_0,\omega_0,0)\right|_{(I-\Pi)\mathbb{R}^2}\right)^{-1}F_\delta(\theta_0,\omega_0,0),
\end{array}
\end{equation}
is invertible on $\Pi\mathbb{R}^2.$ Then, there exists a unique $(\theta_1,\omega_1)\in\mathbb{R}^2$ such that, for all $|\delta|\not=0$ sufficiently small, one can find $(\theta_{1,\delta},\omega_{1,\delta})\in\mathbb{R}^2$ that satisfies both
$$
   F(\theta_0+\delta\theta_{1,\delta},\omega_0+\delta\omega_{1,\delta},\delta)=0,
$$
and
$$
   (\theta_{1,\delta},\omega_{1,\delta})\to(\theta_1,\omega_1)\quad{\rm as}\quad \delta\to 0.
$$
\end{theorem}

\begin{theorem}\label{theoremkmn2} {\rm (two-dimensional version of \cite[Theorem 2]{kmn})} Assume all the conditions of Theorem~\ref{theoremkmn}. Let $(\theta_{1,\delta},\omega_{1,\delta})$ be as given by Theorem~\ref{theoremkmn}. Denote by $\lambda_*\in\mathbb{R}$ the eigenvalue of the linear map 
\begin{equation}\label{linearmap}
{\left.\Pi F_{(\theta_0,\omega_0)^2}(\theta_1,\omega_1)^T\right|}_{\Pi\mathbb{R}^2}+\left.\Pi (F_\delta)_{(\theta,\omega)}(\theta_0,\omega_0,0)\right|_{\Pi\mathbb{R}^2}.
\end{equation} Then
$$
  \lambda_\delta=\delta\lambda_*+o(\delta).
$$
\end{theorem}

\vskip0.2cm

\noindent In order to apply Theorem~\ref{theoremkmn} to the Poincare map $(\theta,\omega)\mapsto P(\theta,\omega,\delta)$, we consider
\begin{equation}\label{F}\def\arraystretch{1.2}
\begin{array}{l}
   F(\theta,\omega,\delta)=P(\theta,\omega,\delta)-(\theta,\omega)^T,\\ 
   \Pi\zeta=z^T\zeta y,
\end{array}
\end{equation}
and notice that (\ref{dd}) implies
\begin{equation}\label{ddd}
    z^TP_{(\theta,\omega)^2}(\xi(s),0)y=0,\quad\mbox{for all}\ s\in\mathbb{R},
\end{equation}
which allows (as we show in the proof of Theorem~\ref{theoremmain}), to ignore all the expressions of Theorem~\ref{theoremkmn} that involve the second derivative.

\vskip0.2cm

\noindent{\bf Proof of Theorem \ref{theoremmain}.} {\it The necessity part.} Here we follow the idea of  Makarenkov-Ortega \cite[Lemma~2]{ortega}.
Assume that $P(\theta_\delta,\omega_\delta,\delta)=(\theta_\delta,\omega_\delta)^T$, $\delta\in\mathbb{R}$, for some family $\left\{(\theta_\delta,\omega_\delta)\right\}_{\delta\in\mathbb{R}}$ satisfying (\ref{convergence}). We claim that (\ref{necessary}) holds.

\vskip0.2cm

\noindent The derivative $F'(\theta,\omega,\delta)$  of the $C^1$ function (\ref{F}) is a $2\times 3$-matrix. Observe that ${\rm rank}\hskip0.05cm F'(\xi(s_0),0)=1.$ Otherwise the equation $F(\theta,\omega,\delta)=0$ should describe a curve in a small neighborhood of $(\xi(s_0),0)$. However, the set $\{(\theta,\omega,\delta):F(\theta,\omega,\delta)=0\}$ contains both the curve $\{(\xi(s),0)\}_{s\in\mathbb{R}}$ and also the set $\{(\theta_\delta,\omega_\delta,\delta)\}_{\delta\in\mathbb{R}}.$ Now we know that 
${\rm rank}\hskip0.05cm F'(\xi(s_0),0)=1$ and it remains to prove that 
\begin{equation}\label{remains}
  {\rm rank}\hskip0.05cm F'(\xi(s_0),0)=2,\quad{\rm if}\quad z^TF_\delta(\xi(s_0),0)\not=0.
\end{equation}
By Fredholm alternative for matrices (see e.g. \cite[Theorem~4.5.3]{kuttler}),
$$
   {\rm Im}\hskip0.05cm F_{(\theta,\omega)}(\xi(s_0),0)=\left({\rm Ker}\hskip0.05cm F_{(\theta,\omega)}(\xi(s_0),0)^T\right)^\perp
$$
Since ${\rm Ker}\hskip0.05cm \Phi_{(\theta,\omega)}(\xi(s_0),0)^T={\rm span}(z)$, we conclude that $\left({\rm Ker}\hskip0.05cm F_{(\theta,\omega)}(\xi(s_0),0)^T\right)^\perp={\rm span}(\tilde y)$, where $\tilde y$ is an eigenvector of  $F_{(\theta,\omega)}(\xi(s_0),0)$ that corresponds to the non-zero eigenvalue of $F_{(\theta,\omega)}(\xi(s_0),0)$. Therefore, ${\rm Im}\hskip0.05cm F_{(\theta,\omega)}(\xi(s_0),0)={\rm span}(\tilde y)$. But $z^T F_\delta(\xi(s_0),0)\not=0$ implies, see formula~(\ref{decomposition}), that the vectors $\tilde y$ and 
$F_\delta(\xi(s_0),0)$ are linearly independent, which completes the proof of (\ref{remains}).

\vskip0.2cm

\noindent{\it The sufficiency part.}  Here we use Theorem~\ref{theoremkmn}. The projector $\Pi$ defined by (\ref{F}) is invariant with respect to $F_{(\theta,\omega)}(\theta_0,\omega_0,0)$ and the projector $I-\Pi$ is given by, see formula (\ref{decomposition}),
$$
   (I-\Pi)\zeta=\tilde z^T\zeta \tilde y,
$$
so that $F_{(\theta,\omega)}(\theta_0,\omega_0,0)$ is invertible on $(I-\Pi)\mathbb{R}^2.$ The requirement $\Pi F_\delta(\theta_0,\omega_0,0)=0$ of Theorem~\ref{theoremkmn} holds by (\ref{necessary}), and the requirement $\Pi F_{(\theta,\omega)^2}(\theta_0,\omega_0,0)\Pi h_1\Pi h_2=0$ holds by (\ref{ddd}). The properties (\ref{F}) and (\ref{ddd}) imply that the expression  (\ref{invertible}) is invertible on ${\rm span}(y)$ if and only if (\ref{sufficient}) holds. Therefore, the conclusion of the theorem follows by applying Theorem~\ref{theoremkmn}.

\vskip0.2cm

\noindent {\it The stability part.} Assume that conditions (\ref{stab1}) and (\ref{stab2}) hold. Let $\rho_\delta$ be the eigenvalue of $P_{(\theta,\omega)}(\theta_\delta,\omega_\delta,\delta)$ such that
$$
   \rho_\delta\to 1\quad{\rm as}\quad \delta\to 0.
$$
We have to show that $|\rho_\delta|<1$ for all $|\delta|>0$ sufficiently small. Observe that $$\lambda_\delta=\rho_\delta-1$$ is the eigenvalue of $F_{(\theta,\omega)}(\theta_\delta,\omega_\delta,\delta)$. As it was established in the sufficiency part of the proof, the expression (\ref{linearmap}) coincides with $z^T (P_\delta)_{(\theta,\omega)}(\theta_0,\omega_0,0)y$. Therefore,  condition (\ref{stab2}) ensures that $\lambda_*$ of Theorem~\ref{theoremkmn2} verifies $\lambda_*<0$ and so 
Theorem~\ref{theoremkmn2} ensures that $\lambda_\delta<0$ for all $\delta>0$ sufficiently small.

\vskip0.2cm

\noindent The proof of the theorem is complete.\qed

\section{The solution of equation (\ref{linearization})}
\label{hfsolution}
 The solution $(h(t),f(t))$ of equation (\ref{linearization}) is given by
\begin{eqnarray*}
h(t)&=& \dfrac{1}{384} e^{-3 t} [384 e^{3 t} + (\omega - \theta)^3 - 
   e^{6 t} (\omega + \theta)^3 + 
   e^{2 t} \{-192 + 3 \omega^3 (3 + 4 t) +\\
   &&+ 
      3 \omega^2 (1 - 4 t) \theta - 
      3 \omega (7 + 4 t) \theta^2 + (1 + 12 t) \theta^3\} + 
   e^{4 t} (-192 + 3 \omega^3 (-3 + 4 t) - \\
   &&-
      3 \omega (-7 + 4 t) \theta^2 + (1 - 12 t) \theta^3 + 
      3 \omega^2 (\theta + 4 t \theta))],
\end{eqnarray*}
\begin{eqnarray*}
f(t)&=& \dfrac{1}{7680}e^{-3 t} [-1920 e^{2 t} - 1920 e^{4 t} - 56 \omega^3 + 
   60 e^{2 t} \omega^3 - 60 e^{4 t} \omega^3 + 56 e^{6 t} \omega^3 +\\
   &&+ 
   120 e^{2 t} \omega^3 t 
   +120 e^{4 t} \omega^3 t + 
   168 \omega^2 \theta + 60 e{2 t} \omega^2 \theta + 
   60 e^{4 t} \omega^2 \theta + 168 e^{6 t} \omega^2 \theta - \\
   && -
   120 e^{2 t} \omega^2 t \theta + 120 e^{4 t} \omega^2 t \theta - 
   168 \omega \theta^2 - 780 e^{2 t} \omega \theta^2 + 
   780 e^{4 t} \omega \theta^2 + 168 e^{6 t} \omega \theta^2 - \\
   &&-
   120 e^{2 t} \omega t \theta^2 - 120 e^{4 t} \omega t \theta^2 + 
   56 \theta^3 + 580 e^{2 t} \theta^3 + 580 e^{4 t} \theta^3 + 
   56 e^{6 t} \theta^3 + \\
   &&+120 e^{2 t} t \theta^3 - 
   120 e^{4 t} t \theta^3 + 3 e^t \{-65 (\omega - 3 \theta) (\omega - \theta)^2 +       65 e^{4 t} (\omega + \theta)^2 (\omega + 3 \theta) +\\
   &&+       e^{2 t} (1280 + 140 \omega^3 t - 921 \omega^2 \theta +          60 \omega t \theta^2 - 697 \theta^3)\} \cos t +  12 e^{2 t} \{(-1 + e^{2 t}) \omega^3 + \\
   &&+      13 (1 + e^{2 t}) \omega^2 \theta +       3 (-1 + e^{2 t}) \omega \theta^2 - 
      9 (1 + e^{2 t}) \theta^3\} \cos(2 t) +    45 e^{3 t} \omega^2 \theta \cos(3t) -\\
      &&-    135 e^{3 t} \theta^3 \cos(3t) - 195 e^t \omega^3 \sin t - 
   1179 e^{3 t} \omega^3 \sin t - 195 e^{5 t} \omega^3 \sin t  -\\
   &&-    195 e^t \omega^2 \theta \sin t +    195 e^{5 t} \omega^2 \theta \sin t +    1260 e^{3 t} \omega^2 t \theta \sin t +    975 e^t \omega \theta^2 \sin t +  \\
   &&+  3813 e^{3 t} \omega \theta^2 \sin t +    975 e^{5 t} \omega \theta^2 \sin t - 585 e^t \theta^3 \sin t  +    585 e^{5 t} \theta^3 \sin t + \\
   &&+540 e^{3 t} t \theta^3 \sin t  -    24 e^{2 t} \omega^3 \sin (2t) - 24 e^{4 t} \omega^3 \sin(2t) -    48 e^{2 t} \omega^2 \theta \sin(2t)  +  \\
   &&+  48 e^{4 t} \omega^2 \theta \sin(2t) + 
   288 e^{2 t} \omega \theta^2 \sin(2t)  +    288 e^{4 t} \omega \theta^2 \sin(2t) -    216 e^{2 t} \theta^3 \sin (2t) + \\
   &&+216 e^{4 t} \theta^3 \sin(2t)  -    5 e^{3t} \omega^3 \sin(3t)  + 135 e^{3 t} \omega \theta^2 \sin(3t)].
\end{eqnarray*}

\section*{Compliance with Ethical Standards}
\noindent {\bf Conflict of Interest:} The authors have no conflict of interest.

\bibliographystyle{siamplain}

\end{document}